\documentclass[a4paper, 11pt]{article}

\usepackage{mathrsfs}
\usepackage{epsfig}
\usepackage{amsmath}
\usepackage{amssymb}
\usepackage{latexsym}
\usepackage{amsfonts}
\usepackage{amsthm}
\usepackage[numbers,sort&compress]{natbib}
\usepackage{graphicx}
\ifx\pdfoutput\undefined  \else
 \fi
\usepackage[centerlast]{caption2}
\usepackage{color}

\newtheorem{definition}{Definition}

\newtheorem{lemma}{Lemma}
\newtheorem{proposition}{Proposition}
\newtheorem{theorem}{Theorem}
\newtheorem{observation}{Observation}

\numberwithin{figure}{section}
\numberwithin{definition}{section}
\numberwithin{observation}{section}
\numberwithin{lemma}{section}
\numberwithin{theorem}{section}

\begin{document}

\title{
{The crossing number of locally twisted cubes} \footnote{The
research is supported by NSFC (60973014, 60803034, 11001035) and
SRFDP (200801081017, 200801411073).}
\author{
\ Haoli Wang, \ Xirong Xu, \ Yuansheng Yang\footnote {corresponding
author's email : yangys@dlut.edu.cn}, \ Bao Liu\\
Department of Computer Science, \\
Dalian University of Technology, Dalian, 116024, P. R. China\\
\\
Wenping Zheng \\
Key Laboratory of Computational Intelligence and Chinese Information\\
Processing of Ministry of Education,\\
Shanxi University, Taiyuan, 030006, P. R. China\\
\\
Guoqing Wang \\
Center for Combinatorics, LPMC-TJKLC \\
Nankai University, Tianjin, 300071, P. R. China\\
}}

\date{}
\maketitle
\begin{abstract}

The {\it crossing number} of a graph $G$ is the minimum number of
pairwise intersections of edges in a drawing of $G$. Motivated by
the recent work [Faria, L., Figueiredo, C.M.H. de, Sykora, O.,
Vrt'o, I.: An improved upper bound on the crossing number of the
hypercube. J. Graph Theory {\bf 59}, 145--161 (2008)] which solves
the upper bound conjecture on the crossing number of $n$-dimensional
hypercube proposed by Erd\H{o}s and Guy, we give upper and lower
bounds of the crossing number of locally twisted cube, which is one
of variants of hypercube.

\end{abstract}

\section{Introduction}

\indent \indent  The {\it crossing number} $cr(G)$ of a graph $G$ is
the minimum number of pairwise intersections of edges in a drawing
of $G$ in the plane. The notion of crossing number is a central one
for Topological Graph Theory and has been studied extensively  by
mathematicians including Erd\H{o}s, Guy, Tur\'{a}n and Tutte, et al.
(see
\cite{EG73,Guy60,PaSpTo00,PachToth00,SSSV97,WB78,Turan77,Tutte70}).
In the past thirty years, it turned out that crossing number has
many important applications in discrete and computational geometry
(see \cite{Bi91,M02,PaSh98,Sz97,SoTaTo02,SoTo01,TaoVu}). For
example, Sz\'{e}kely \cite{Sz97} employed the `crossing lemma'
\cite{ACNS82,L83}
 to give a simple proof of the following
well-known theorem in discrete and computational geometry.

{\bf Theorem A.} (Szemer\'{e}di-Trotter \cite{SzemTro83}) \ {\sl
Given n points and $\ell$ lines in the plane, there is a constant
$c$ for which the number of incidences among the points and lines is
at most $c[(n\ell)^{2/3}+n+\ell]$.}

On the other hand, the immediate applications in VLSI theory and
wiring layout problems (see \cite{BL84,S05,L81,L83}) also inspired
the study of crossing number of some popular parallel network
topologies such as hypercube and its variations. Among all the
popular parallel network topologies, hypercube is the first to be
studied (see \cite{DR95,Egg70,FF00,FFSV08,M91,SV93}). An
$n$-dimensional hypercube $Q_n$ is a graph in which the nodes can be
one-to-one labeled with 0-1 binary sequences of length $n$, so that
the labels of any two adjacent nodes differ in exactly one bit.

Computing the crossing number was proved to be NP-complete by Garey
and Johnson \cite{GJ83}. Thus, it is not surprising that the exact
crossing numbers are known for graphs of few families and that the
arguments often strongly depend on their structures (see for example
\cite{Fiorini86,PanRich07,RichterThomassen95,Zarankiewicz54}). Even
for hypercube, for a long time the only known result on the exact
value of crossing number of $Q_n$ has been $cr(Q_3)=0$, $cr(Q_4)=8$
\cite{DR95}, $cr(Q_5)\leq 56$ \cite{M91}. Hence, it is more
practical to find upper and lower bounds of crossing numbers of some
kind of graphs. Concerned with upper bound of crossing number of
hypercube, Erd\H{o}s and Guy \cite{EG73} in 1973 conjectured the
following:
$$cr(Q_n)\leq \frac{5}{32}4^n-\lfloor\frac{n^2+1}{2}\rfloor
2^{n-2}.$$ In 2008, Faria, Figueiredo, Sykora and Vrt'o
\cite{FFSV08} constructed a drawing of $Q_n$ in the plane which has
the conjectured number of crossings mentioned above. Early in 1993
Sykora and Vrt'o \cite{SV93} also proved a lower bound of $cr(Q_n)$:
$$cr(Q_n)>\frac{1}{20}4^n-(n^2+1)2^{n-1}.$$

Since the hypercube does not have the smallest possible diameter for
its resources, to achieve smaller diameter with the same number of
nodes and links as an $n$-dimensional cube, a variety of hypercube
variants were proposed. Locally twisted cube is one of these
variants. The {\it $n$-dimensional locally twisted cube} $LTQ_n$,
proposed by Yang et al. \cite{YEM05} in 2005, keeps as many nice
properties of hypercube as possible and is conceptually closer to
traditional hypercube, while it has diameters of about half of that
of a hypercube of the same size. Therefore, it would be more
attractive to study the crossing number of the $n$-dimensional
locally twisted cubes.

The {\it $n$-dimensional locally twisted cube} $LTQ_n(n\geq 2)$ is defined recursively as follows. \\
\indent (a) $LTQ_2$ is a graph isomorphic to $Q_2$. \\
\indent (b) For $n\geq 3$, $LTQ_n$ is built from two disjoint copies
of $LTQ_{n-1}$ according to the following steps. Let $0LTQ_{n-1}$
denote the graph obtained by prefixing the label of each vertex of
one copy of $LTQ_{n-1}$ with 0, let $1LTQ_{n-1}$ denote the graph
obtained by prefixing the label of each vertex of the other copy
$LTQ_{n-1}$ with 1, and connect each vertex $x=0x_2x_3\ldots x_n$ of
$0LTQ_{n-1}$ with the vertex $1(x_2+x_n)x_3\ldots x_n$ of
$1LTQ_{n-1}$ by an edge, where $‘+’$ represents the modulo 2
addition.

The graphs shown in Figure 1.1 are $LTQ_3$ and $LTQ_4$,
respectively.
\begin{figure}[ht]
\centering
\includegraphics[scale=0.9]{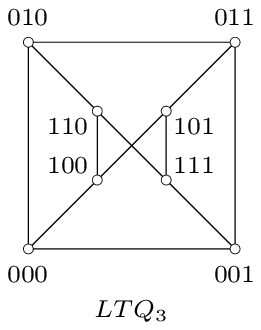} \hspace{40pt}
\includegraphics[scale=0.9]{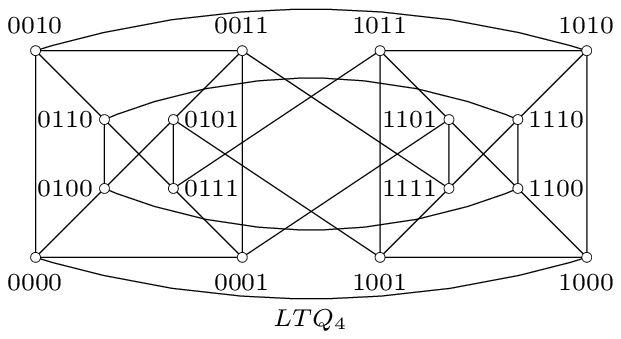}
\caption{\small{Locally twisted cubes $LTQ_3$ and $LTQ_4$}}
\end{figure}

In this paper, we mainly obtain the following bounds of the crossing number
of $LTQ_n$:

$$\frac{4^n}{20}-(n^2+1)2^{n-1}<cr(LTQ_n)\leq\frac{265}{6}4^{n-4}-(n^2+\frac{15+(-1)^{n-1}}{6})2^{n-3}.$$

\section{Upper bound for $cr(LTQ_n)$}

\indent \indent A drawing of $G$ is said to be a {\it good} drawing,
provided that no edge crosses itself, no adjacent edges cross each
other, no two edges cross more than once, and no three edges cross
in a point. It is well known that the crossing number of a graph is
attained only in {\it good} drawings of the graph. So, we always
assume that all drawings throughout this paper are good drawings.
For a good drawing $D$ of a graph $G$, let $\nu_D(G)$ be the number
of crossings in $D$. In what follows, $\nu_D(G)$ is abbreviated to
$\nu_D$ when it is unambiguous.

Let $x=x_1x_2\cdots x_n$ and $y=y_1y_2\cdots y_n$ be two vertices of
$LTQ_n$. Denote
$$\mathscr{D}(x_1x_2\cdots x_n)=2^{n-1}x_1+2^{n-2}x_2+\cdots
+2^0x_n$$ to be the corresponding decimal number of $x_1x_2\cdots
x_n$. Let
$$\theta_i(x)=x_i\ \ \  \mbox{ for }  i\in\{1,2,\ldots,n\}.$$
Let $\lambda(x,y)$ be the smallest positive integer $i\in
\{1,2,\ldots,n\}$ such that $\theta_i(x)\neq \theta_i(y)$. We define
$$\begin{array}{llll}Dim(x,y)=\left\{\begin{array}{llll}
                \lambda(x,y),  & \mbox{if \ \ } x\mbox{ and }y \mbox{ are adjacent};\\
                \infty,  & \mbox{otherwise. \ \ } \\
              \end{array}
              \right .
\end{array}$$
In particular, for an edge $e=xy$, let $Dim(e)=Dim(x,y)$ and say the
edge $e$ lies in the $Dim(e)$-dimension. We call $x$ an {\it odd
vertex} if $|\{1\leq i\leq n: x_i=1\}|\equiv 1\pmod 2$, and an {\it
even vertex} if otherwise.

For the clearness of composition, in the rest of this section,  any
vertex $x\in V(LTQ_n)$ in figures will be represented by the
corresponding decimal number $\mathscr{D}(x)$. We first give a
drawing of $LTQ_4$ with 10 crossings and a drawing of $LTQ_5$ with
68 crossings as shown in Figure 2.1. Hence, we have the following

\begin{proposition}\label{Theorem cr(LTQ5) leq 68}
$cr(LTQ_4)\leq 10$ and $cr(LTQ_5)\leq 68$.
\end{proposition}
\begin{figure}[ht]
\centering
\includegraphics[scale=0.9]{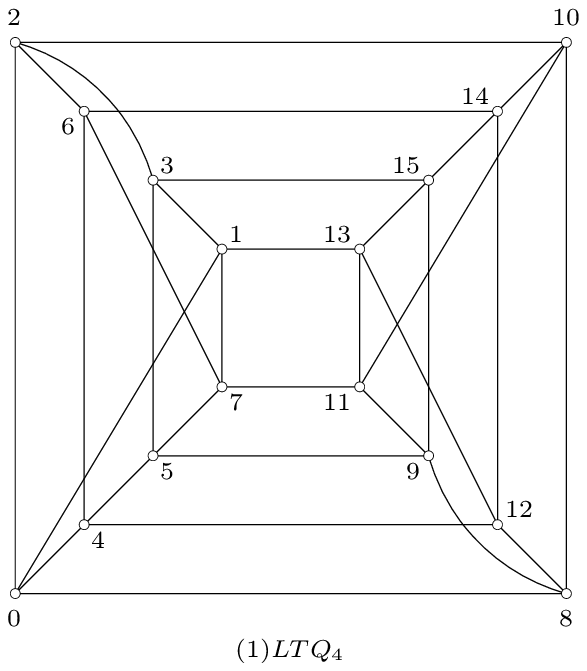} \hspace{40pt}
\includegraphics[scale=0.9]{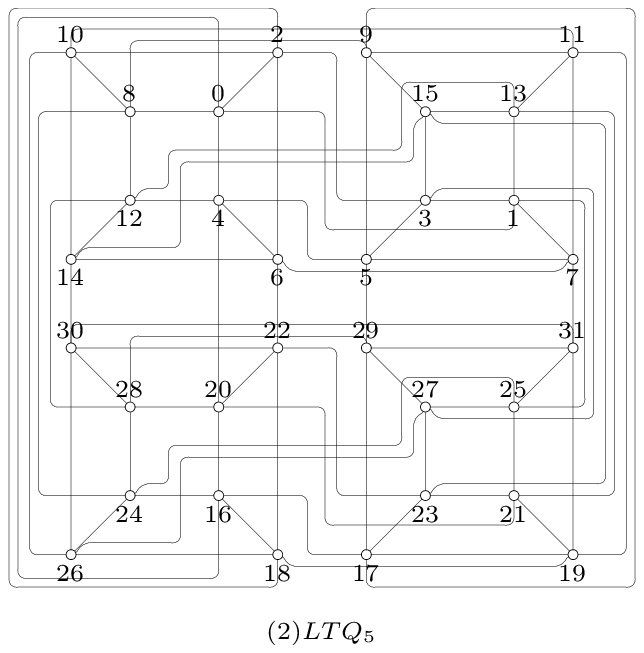}
\caption{\small{Drawings of $LTQ_4$ with 10 crossings and $LTQ_5$
with 68 crossings}}
\end{figure}

Before proving the upper bound of $cr(LTQ_n)$ for $n\geq 6$, we need
to introduce some technical notations. We define two structures
$M^i$ and $M_c^i$, called $``meshes"$ which will be used in counting
the number of crossings. Consider the canonical geometry of the real
plane $\mathbb{R}^2$. By $[0,1]$ we denote the closed interval
joining the points $(0,0)$ and $(1,0)$ of the horizontal real axis.
Let $r$ and $s$ be a non-horizontal pair of parallel straight lines
in the real plane $\mathbb{R}^2$, such that the point $(0,0)$
belongs to $r$ and the point $(1,0)$ belongs to $s$. For a positive
integer $n$, let $\mathscr{L}_n=\{(r_i,s_i):i\in\{1,2,\ldots,n\}\}$
be a set of non-horizontal pairs of parallel straight lines in the
real plane $\mathbb{R}^2$, such that the point $(0,0)$ belongs to
$r_i$ and the point $(1,0)$ belongs to $s_i$.

A {\it mesh} with index $n$, denoted $M^n$, is the set of points of
the plane consisting of the points of the $n$-element set
$\mathscr{L}_n$ plus the points in the interval $[0,1]$. In Figure
2.2, we show as an example a drawing of each $M^1$, $M^2$, $M^3$ and
$M^5$.
\begin{figure}[ht]
\centering
\includegraphics[scale=1.0]{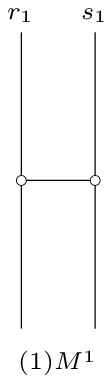} \hspace{10pt}
\includegraphics[scale=1.0]{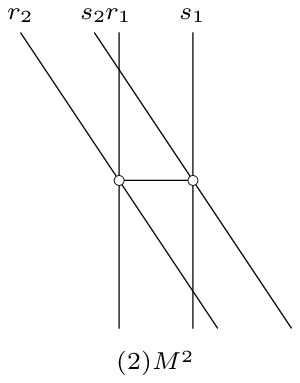} \hspace{10pt}
\includegraphics[scale=1.0]{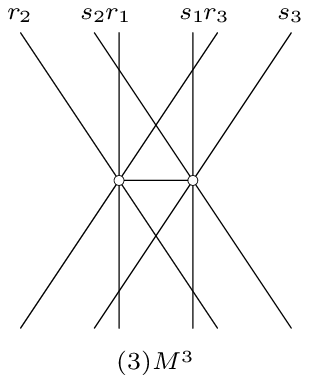} \hspace{10pt}
\includegraphics[scale=1.0]{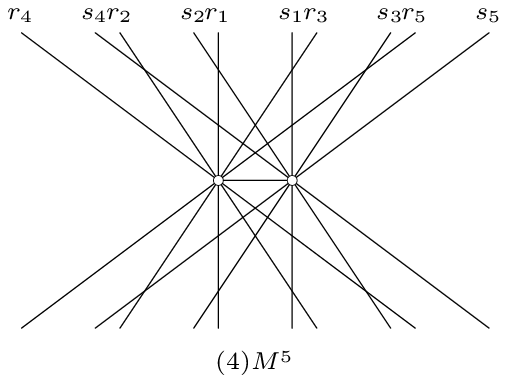}
\caption{\small{Drawings of $M^1$, $M^2$, $M^3$ and $M^5$}}
\end{figure}

A {\it chopped mesh} with index $n$, denoted $M_c^n$, is the set of
points of $M^n$ without a pair of parallel semi-straight lines of
the left-most lower semi-plane. In Figure 2.3, we show a drawing of
each $M_c^1$, $M_c^2$, $M_c^3$ and $M_c^5$.
\begin{figure}[ht] \label{Figure Mc}
\centering
\includegraphics[scale=1.0]{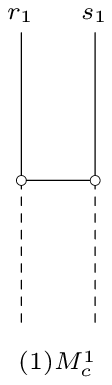} \hspace{10pt}
\includegraphics[scale=1.0]{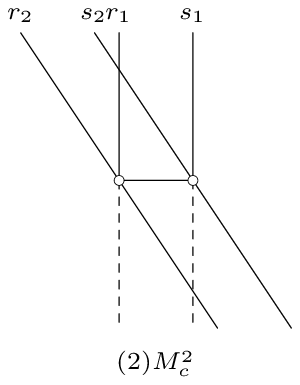} \hspace{10pt}
\includegraphics[scale=1.0]{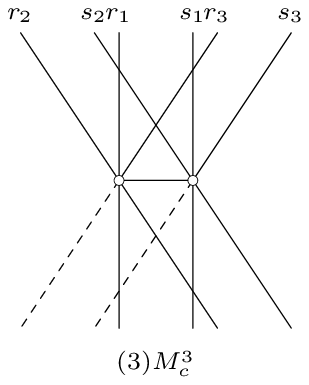} \hspace{10pt}
\includegraphics[scale=1.0]{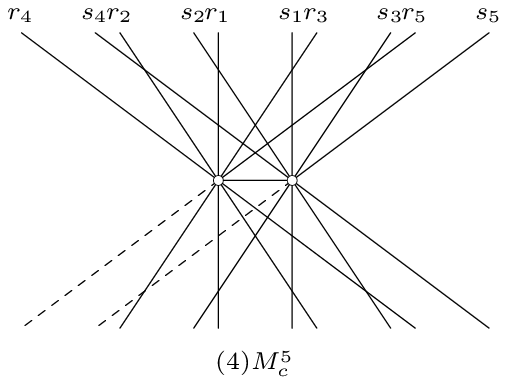}
\caption{\small{Drawings of $M_c^1$, $M_c^2$, $M_c^3$ and $M_c^5$}}
\end{figure}

\begin{lemma}\cite{FFSV08}\label{Lemma mesh} For any positive
integer $n$, there is a drawing of $M^n$ with $n(n-1)$ crossings.
\end{lemma}

\begin{lemma}\cite{FFSV08}\label{Lemma mesh c} For any positive
integer $n$, there is a drawing of $M_c^n$ with $(n-1)^2$ crossings.
\end{lemma}

To prove the general upper bound of $cr(LTQ_n)$, we need to
construct a drawing $D_n$ of $LTQ_n$ with the desired number of
crossings. The philosophy is putting the obtained drawing $D_{n-1}$
of $LTQ_{n-1}$ on the given coordinate systems (see Figure 2.5) and
then replacing each vertex of $LTQ_{n-1}$ by two vertices of $LTQ_n$
and replacing each edge of $LTQ_{n-1}$ by a bunch of two edges of
$LTQ_n$. Hence, we need the following definitions.

\begin{definition} Let $x$ be a vertex of $LTQ_n$, and let $e\in E(LTQ_n)$ be an edge incident with $x$.
Assume that $x$ is drawn precisely on some axis $A$. We call $e$ an
\texttt{a-arc} or \texttt{b-arc} with respect to $x$, provided that
the edge $e$ is drawn to be upward from $A$ (based upon the positive
direction of the axis $A$) or to be downward from $A$, respectively.
In particular, let $$\alpha(x)=|\{e\in E(LTQ_n):e \mbox{ is an {\it
a-arc} with respect to }x\}|$$ and
$$\beta(x)=|\{e\in E(LTQ_n):e \mbox{ is a {\it b-arc} with respect to
}x\}|.$$
\end{definition}

For example, as shown in Figure 2.5, the three edges joining vertex
23 and vertices 17, 27, 21 are a-arcs with respect to vertex 23, and
the three edges joining vertex 23 and vertices 22, 39, 15 are b-arcs
with respect to vertex 23.

\begin{definition} Let $x$ and $y$ be two vertices of $LTQ_n$ with $Dim(x,y)=n-1$. Assume that
 $x$ and $y$ are drawn next to each other on some axis. Then
we define the the forward direction of $x$ to be coincident with the
direction from $y$ to $x$ if $\theta_n(x)=1$ and $x$ is an odd
vertex and that the forward direction of $x$ to be coincident with
the direction from $x$ to $y$ if otherwise (see Figure 2.4).
\end{definition}
\begin{figure}[ht]
\centering
\includegraphics[scale=1.0]{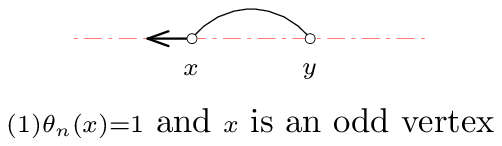} \hspace{80pt}
\includegraphics[scale=1.0]{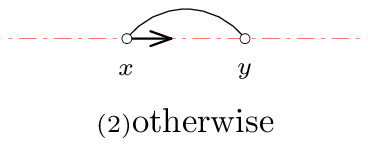}
\caption{\small{The forward direction of vertex $x$}}
\end{figure}

\begin{definition} Let $x$ and $y$ be two adjacent vertices of
$LTQ_n$. For $i\in \{1,2\}$, we define
$\varepsilon_i=\varepsilon_i(x,y)$ and $\zeta_i=\zeta_i(x,y)$
satisfying that
$\{(\varepsilon_1,\zeta_1),(\varepsilon_2,\zeta_2)\}=\{(0,1),(1,0)\}$
if $Dim(x,y)=n-1$ and $\theta_n(x)=1$, and that
$\{(\varepsilon_1,\zeta_1),(\varepsilon_2,\zeta_2)\}=\{(0,0),(1,1)\}$
otherwise.
\end{definition}

In what follows, $\varepsilon_i(x,y),\zeta_i(x,y)$ are abbreviated
to $\varepsilon_i,\zeta_i$ respectively when it is unambiguous. Let
$x=x_1x_2\cdots x_n$ be a vertex of $LTQ_n$. We define
$$x^{\delta}=x_1x_2\cdots x_{n-1} \delta
 x_n$$ to be a vertex of $LTQ_{n+1}$, where $\delta\in \{0,1\}$.

\begin{observation}\label{Observation enlarging} Let $x$ and $y$ be two adjacent vertices of
$LTQ_n$. Then $x^{\varepsilon_i}$ and $y^{\zeta_i}$ are adjacent
vertices of $LTQ_{n+1}$, in particular,
$$\begin{array}{llll}Dim(x^{\varepsilon_i},y^{\zeta_i})=\left\{\begin{array}{llll}
                Dim(x,y),  & \mbox{if \ \ } Dim(x,y)\leq n-1;\\
                n+1,  & \mbox{if \ \ } Dim(x,y)=n;\\
              \end{array}
              \right .
\end{array}$$
\end{observation}

\begin{observation}\label{Observation for the adjacence} Let $x,y,u,v$ be four vertices
of $LTQ_n$ with $Dim(x,u)=Dim(y,v)=n-1$. If $x$ and $y$ are
adjacent, then $u$ and $v$ are adjacent, in particular,
$Dim(u,v)=Dim(x,y)$.
\end{observation}

Now we are in a position to prove the general upper bound of
$cr(LTQ_n)$.

\begin{theorem}\label{Theorem Upper Bound for general n} For
$n\geq 6$,
$$cr(LTQ_n)\leq\frac{265}{6}4^{n-4}-(n^2+\frac{15+(-1)^{n-1}}{6})2^{n-3}.$$
\end{theorem}

\begin{proof} To prove the theorem, we shall construct a drawing $D_n$ of
$LTQ_n$ for any $n\geq 6$, which satisfies the following five
properties.

\textbf{Property 1:} {\sl
$\nu_{D_n}=\frac{265}{6}4^{n-4}-(n^2+\frac{15+(-1)^{n-1}}{6})2^{n-3}.$}

\textbf{Property 2:} {\sl Every vertex $x$ of $LTQ_n$ is drawn
precisely on some axis, and moreover, $|\alpha(x)-\beta(x)|\leq 1$.}

\textbf{Property 3:} {\sl Let $x,u$ be two vertices of $LTQ_n$ with
$Dim(x,u)=n-1$. Then $x$ and $u$ are drawn next to each other on the
same axis. Moreover, $\alpha(x)=\alpha(u)$ and $\beta(x)=\beta(u)$.}

\textbf{Property 4:} {\sl Let $x,y,u,v$ be four vertices of $LTQ_n$
with $Dim(x,u)=Dim(y,v)=n-1$. Assume that $x$ and $y$ are adjacent.
Then $xy$ is an {\it a-arc} ({\it b-arc}) with respect to $x$ if and
only if $uv$ is an {\it a-arc} ({\it b-arc}) with respect to $u$. }

\textbf{Property 5:} {\sl Let $x,y,u,v$ be four vertices of $LTQ_n$
with $Dim(x,u)=Dim(y,v)=n-1$. If $Dim(x,y)<n$ then
$\nu_{D_n}(xy,uv)=0$.}

Assume first $n=6$. The drawing $D_6$ is given in Figure 2.5. It is
not hard to check that Properties 2, 3, 4 and 5 hold for $D_6$. We
verify that the number of crossings is
$400=\frac{265}{6}\cdot4^{6-4}-(6^2+\frac{15+(-1)^{6-1}}{6})\cdot2^{6-3}$,
and so Property 1 holds for $D_6$.
\begin{figure}[ht]
\centering \hspace{-18pt}
\includegraphics[scale=1.0]{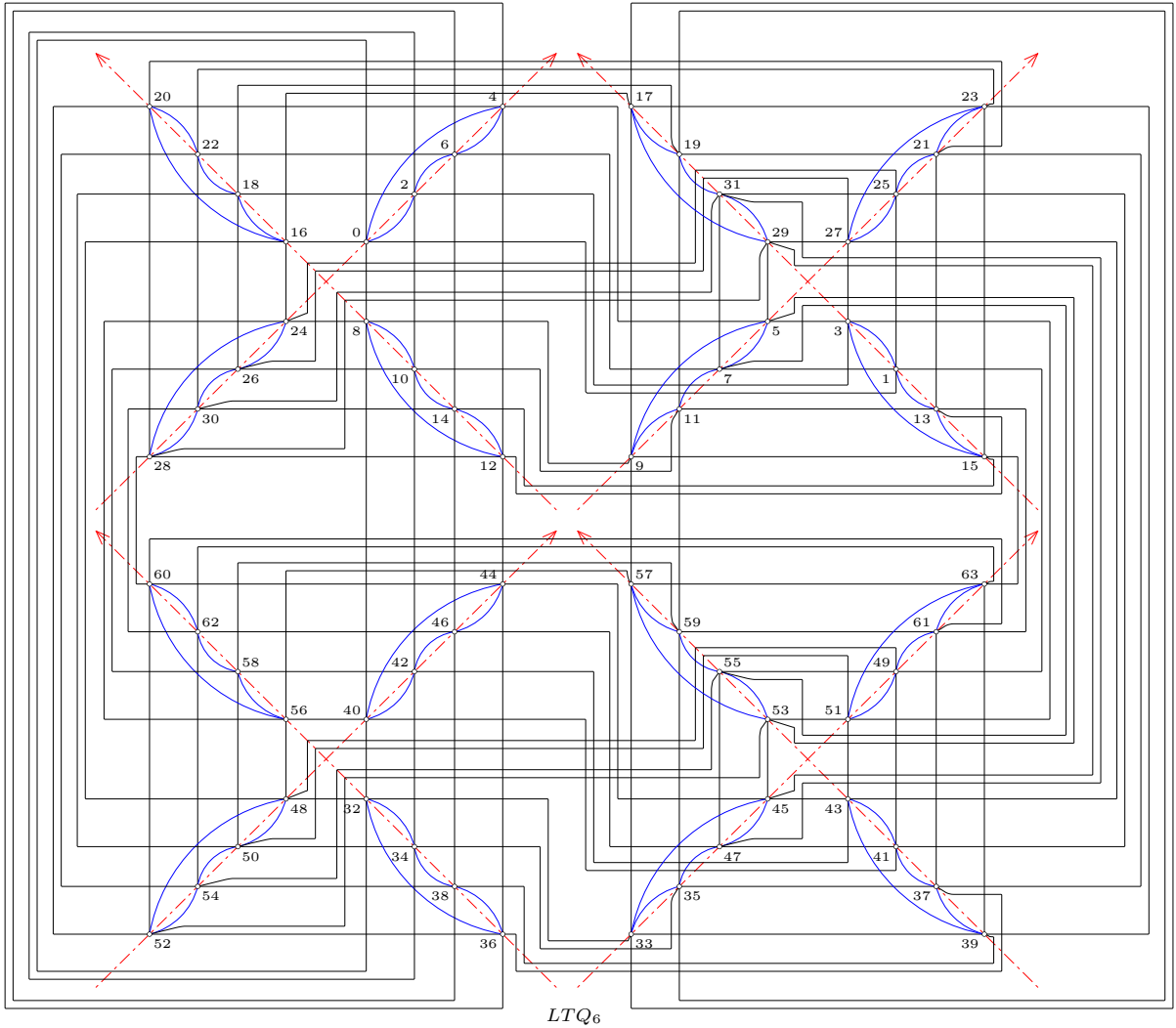}
\caption{\small{The drawing of $D_6$}}
\end{figure}

Now assume that $n\geq 6$ and that there exists a drawing $D_n$ of
$LTQ_n$ satisfying Properties 1, 2, 3, 4 and 5. It suffices to
construct a drawing $D_{n+1}$ of $LTQ_{n+1}$ for which the above
properties hold. The process of constructing $D_{n+1}$ is as
follows. Replace each vertex $x$ of $LTQ_n$ in the ``small"
neighborhood of $x$ in the drawing $D_n$ by two vertices $x^0,x^1\in
V(LTQ_{n+1})$, both of which are drawn precisely on the same axis as
$x$ such that the direction from $x^0$ to $x^1$ is coincident with
the forward direction of $x$. Then join $x^0$ and $x^1$ by an {\sl
a-arc} or {\sl b-arc} with respect to $x^0$ ($x^1$) according to
$\alpha(x)\leq \beta(x)$ or not. By Observation \ref{Observation
enlarging}, we need to replace each edge incident with $x$ in
$LTQ_n$, denoted $e=xy\in E(LTQ_n)$,  by a bunch of two edges
$x^{\varepsilon_1}y^{\zeta_1}, x^{\varepsilon_2}y^{\zeta_2}\in
E(LTQ_{n+1})$ which are ``parallel'' or crossed each other at
``infinity'' (compared to the ``small" neighborhoods of $x$ and
$y$), and drawn along the original edge $e$.

To illustrate the process above, we give in Figure 2.6 the extracted
local drawing on vertices 9, 11, 7, 5 in $D_6$ and the corresponding
extended drawings in $D_7$ and $D_8$. Notice that in Figure 2.6(1)
the vertices 11 and 7 are odd vertices, and that
$Dim(9,11)=Dim(5,7)=5=n-1$. Hence, the forward direction of the
vertex 11(7) is from 9(5) to 11(7).

\begin{figure}[ht]
\centering
\includegraphics[scale=1.0]{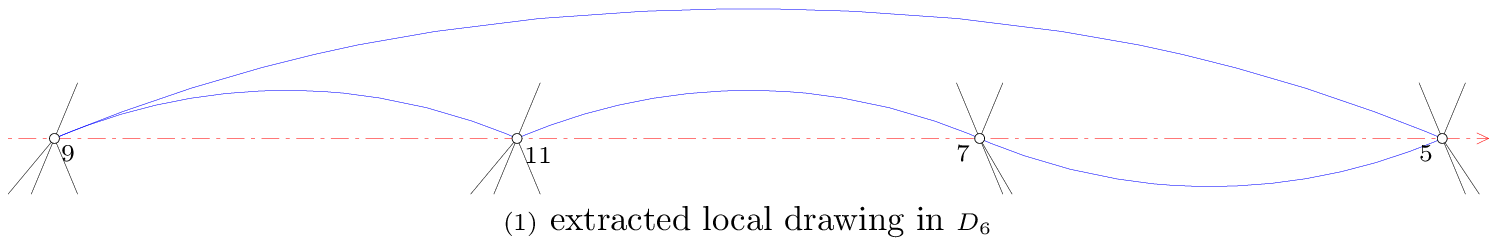}

\vspace{10pt}
\includegraphics[scale=1.0]{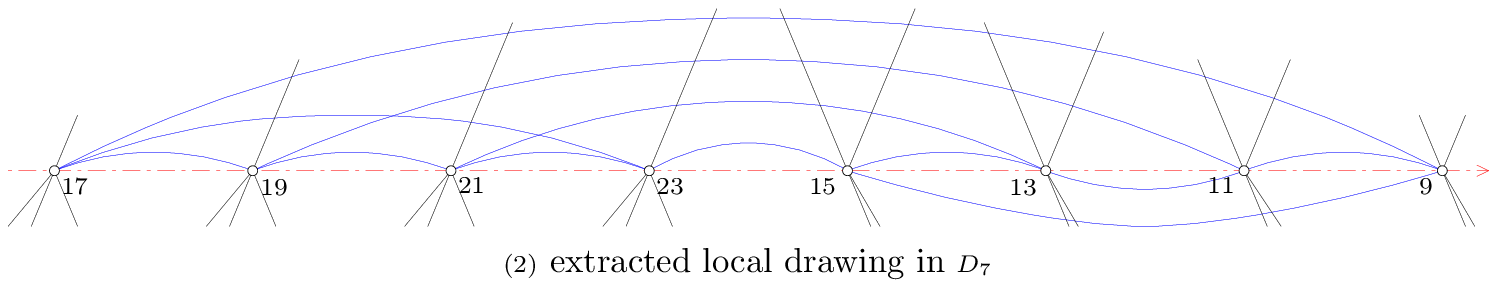}

\vspace{10pt}
\includegraphics[scale=1.0]{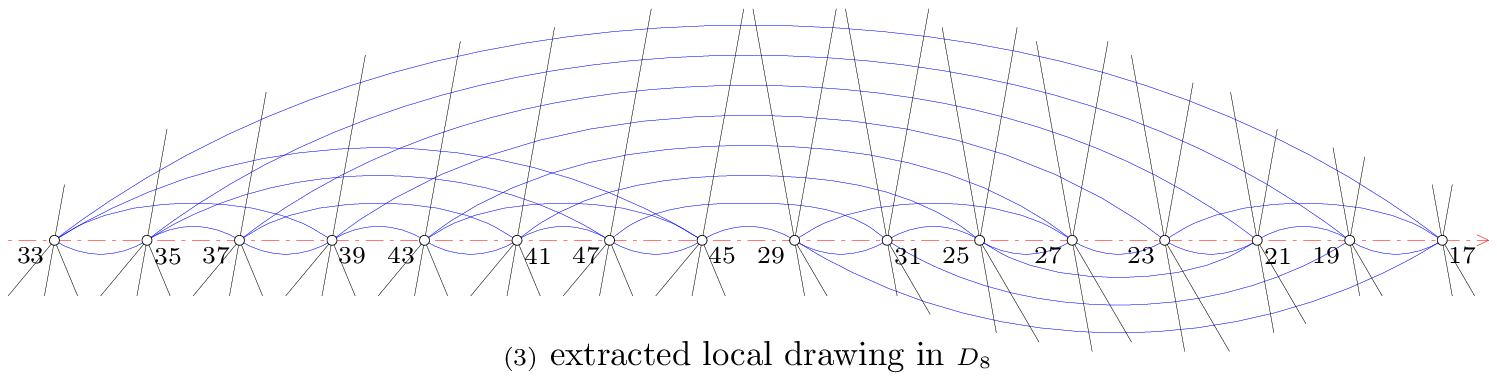}
\caption{\small{The extracted local drawings}}
\end{figure}

By the process described as above, we conclude that Properties 2, 3
and 4 hold for $D_{n+1}$. Because that $D_n$ has Properties 3, 4 and
5, we can verify that
$\nu_{D_{n+1}}(x^{\varepsilon_1}y^{\zeta_1},x^{\varepsilon_2}y^{\zeta_2})=0$
for any edge $xy\in LTQ_n$ with $Dim(xy)<n-1$ (see Figure 2.7, where
$u,v\in V(LTQ_n)$ such that $Dim(u,x)=Dim(v,y)=n-1$) and that
$\nu_{D_{n+1}}(x^{\varepsilon_1}y^{\zeta_1},x^{\varepsilon_2}y^{\zeta_2})=0$
for any edge $xy\in LTQ_n$ with $Dim(xy)=n-1$ (see Figure 2.8).
Combining with Observation \ref{Observation enlarging}, we conclude
that Property 5 holds for $D_{n+1}$.

\begin{figure}[ht]
\centering
\includegraphics[scale=1.0]{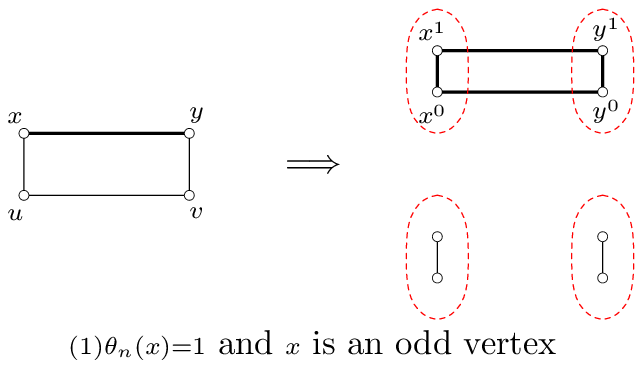} \hspace{30pt}
\includegraphics[scale=1.0]{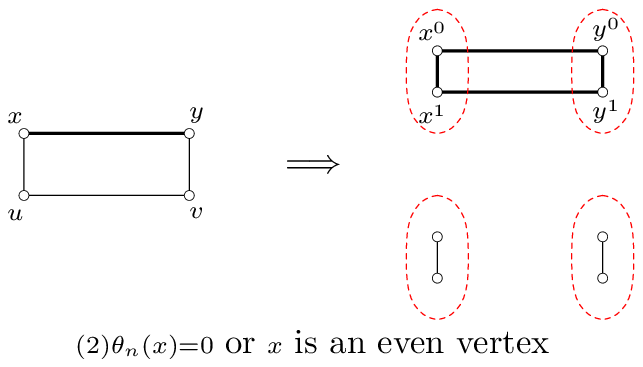}
\caption{\small{The case for $Dim(xy)<n-1$}}
\end{figure}

\begin{figure}[ht]
\centering
\includegraphics[scale=1.0]{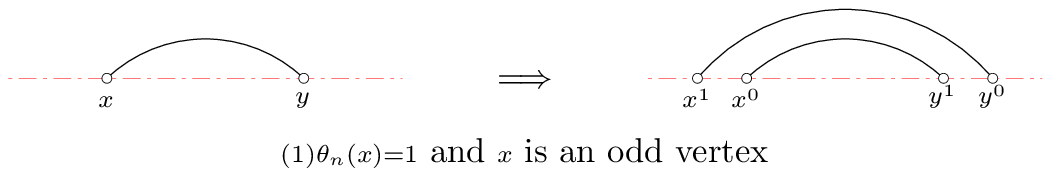}

\vspace{10pt}
\includegraphics[scale=1.0]{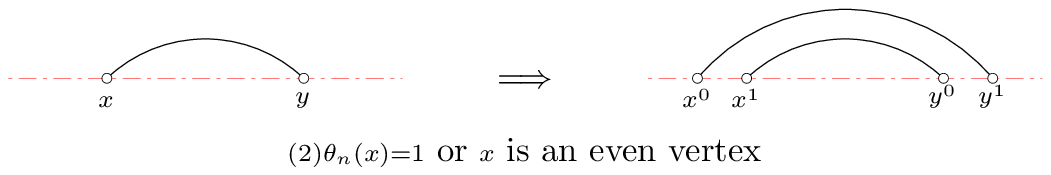}

\vspace{10pt}
\includegraphics[scale=1.0]{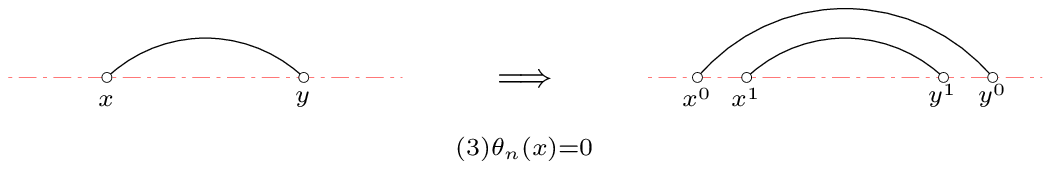}
\caption{\small{The case for $Dim(xy)=n-1$}}
\end{figure}

It remains to show that Property 1 holds for $D_{n+1}$.

\textbf{Claim A.} {\sl For any vertex $x$ of $LTQ_n$, the number of
crossings produced in the ``small" neighborhood of the new edge
$x^{0}x^1$ in $D_{n+1}$ are equal to $\frac{(n-1)^2}{4}$ for odd $n$
and $\frac{n(n-2)}{4}$ for even $n$.}

{\sl Proof of Claim A.} Since $D_{n+1}$ has Properties 2, 3 and 4,
we conclude that the neighborhood of the new edge $x^0x^1$
corresponds to a drawing of $M^{\frac{n+1}{2}}_c$ for odd $n$, and a
drawing of $M^{\frac{n}{2}}$ for even $n$. Then the claim follows
from Lemma \ref{Lemma mesh} and Lemma \ref{Lemma mesh c}. \qed

\textbf{Claim B.} {\sl $|\{xy\in E(LTQ_n): Dim (xy)=n\mbox{ and }
\nu_{D_{n+1}}(x^{\varepsilon_1}y^{\zeta_1},x^{\varepsilon_2}y^{\zeta_2})=1\}|=2^{n-2}.$}

{\sl Proof of Claim B.}  By Observation \ref{Observation for the
adjacence}, there exists a partition $E_1,\ldots, E_{2^{n-2}}$ of
$\{e\in LTQ_n: Dim(e)=n\}$ with $|E_i|=2$, say
$$E_i=\{x_iy_i,u_iv_i\},$$ such that
$$Dim(x_i,u_i)=Dim(y_i,v_i)=n-1,$$
where $i\in \{1,2,\ldots,2^{n-2}\}$. To prove Claim B, it suffices
to show that
\begin{equation}\label{equation crossings equal to 1}
\nu_{D_{n+1}}(x_i^{\varepsilon_1}y_i^{\zeta_1},x_i^{\varepsilon_2}y_i^{\zeta_2})+
\nu_{D_{n+1}}(u_i^{\varepsilon_1}v_i^{\zeta_1},u_i^{\varepsilon_2}v_i^{\zeta_2})=1
\end{equation}
for all $i\in \{1,2,\ldots,2^{n-2}\}$. Assume without loss of
generality that $\theta_n(y_i)=\theta_n(v_i)=1$ and $v_i$ is an odd
vertex, i.e., $\theta_n(x_i)=\theta_n(u_i)=0$ and $y_i$ is an even
vertex. Since $D_n$ has properties 3 and 4, we can verify
\eqref{equation crossings equal to 1} immediately by two cases
$\nu_{D_n}(x_iy_i,u_iv_i)=0$ and $\nu_{D_n}(x_iy_i,u_iv_i)=1$, which
are shown in Figure 2.9. This proves Claim B. \qed
\begin{figure}[ht]
\centering
\includegraphics[scale=1.0]{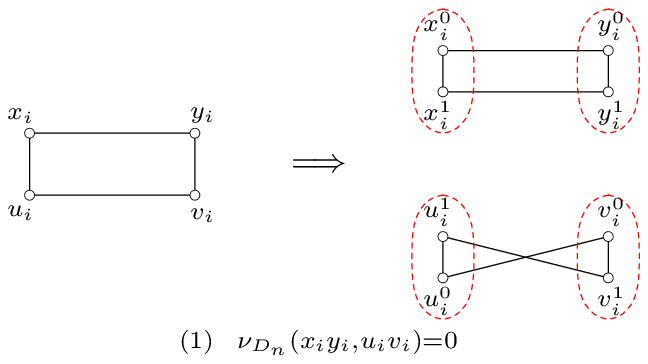} \hspace{30pt}
\includegraphics[scale=1.0]{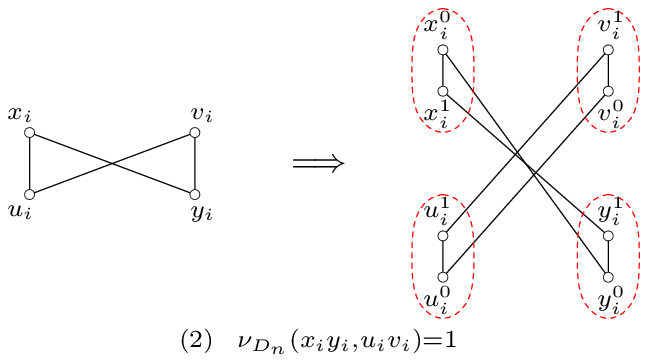}
\caption{\small{Two cases of $\nu_{D_n}(x_iy_i,u_iv_i)=0$ and
$\nu_{D_n}(x_iy_i,u_iv_i)=1$}}
\end{figure}

By the process of constructing $D_{n+1}$, we conclude that
\begin{equation}\label{equation total crossing number}\nu_{D_{n+1}}=4\cdot\nu_{D_{n}}+\Gamma_n+|\{xy\in E(LTQ_n):
\nu_{D_{n+1}}(x^{\varepsilon_1}y^{\zeta_1},x^{\varepsilon_2}y^{\zeta_2})=1\}|
\end{equation}
where $\Gamma_n$ denotes the total number of crossings produced in
the ``small" neighborhoods of all new edges $x^0x^1$. By Claim A, we
have that
\begin{equation}\label{equation number of crossings M}
\begin{array}{llll}\Gamma_n=\left\{\begin{array}{llll}
               2^n\cdot\frac{(n-1)^2}{4},  & \mbox{if \ \ } n\equiv 1 \pmod 2;\\
               2^n\cdot\frac{n(n-2)}{4},  & \mbox{if \ \ } n\equiv 0 \pmod 2.\\
              \end{array}
              \right.
\end{array}
\end{equation}
Recall that $D_{n+1}$ has Property 5. It follows from Observation
\ref{Observation enlarging} that $|\{xy\in E(LTQ_n): Dim (xy)\leq
n-1\mbox{ and }
\nu_{D_{n+1}}(x^{\varepsilon_1}y^{\zeta_1},x^{\varepsilon_2}y^{\zeta_2})=1\}|=0$.
By Claim B, we have that
\begin{equation}\label{equation crossings at infinity}
|\{xy\in E(LTQ_n):
\nu_{D_{n+1}}(x^{\varepsilon_1}y^{\zeta_1},x^{\varepsilon_2}y^{\zeta_2})=1\}|
=2^{n-2}.
\end{equation}

By \eqref{equation total crossing number}, \eqref{equation number of
crossings M} and \eqref{equation crossings at infinity}, we conclude
that
$$\begin{array}{llll}\nu_{D_{n+1}} =\left\{\begin{array}{llll}
               4\cdot\nu_{D_n}+2^n\cdot\frac{(n-1)^2}{4}+2^{n-2}=4\nu_{D_n}+(n^2-2n+2)2^{n-2},  & \mbox{if \ \ } n\equiv 1 \pmod 2;\\
               4\cdot\nu_{D_n}+2^n\cdot\frac{n(n-2)}{4}+2^{n-2}=4\nu_{D_n}+(n^2-2n+1)2^{n-2},  & \mbox{if \ \ } n\equiv 0 \pmod 2.\\
              \end{array}
              \right.
\end{array}$$
Since $D_n$ has Property 1, it is easy to verify that Property 1
holds for $D_{n+1}$. This completes the proof of Theorem
\ref{Theorem Upper Bound for general n}.
\end{proof}

\bigskip

For the convenience of the reader, we offer in Figure 2.10 and 2.11
drawings for $LTQ_7$ and $LTQ_8$ obtained according to the process
of constructing $D_n$.

\newpage

\begin{figure}[ht]
\hspace{-35pt}
\includegraphics[scale=0.8]{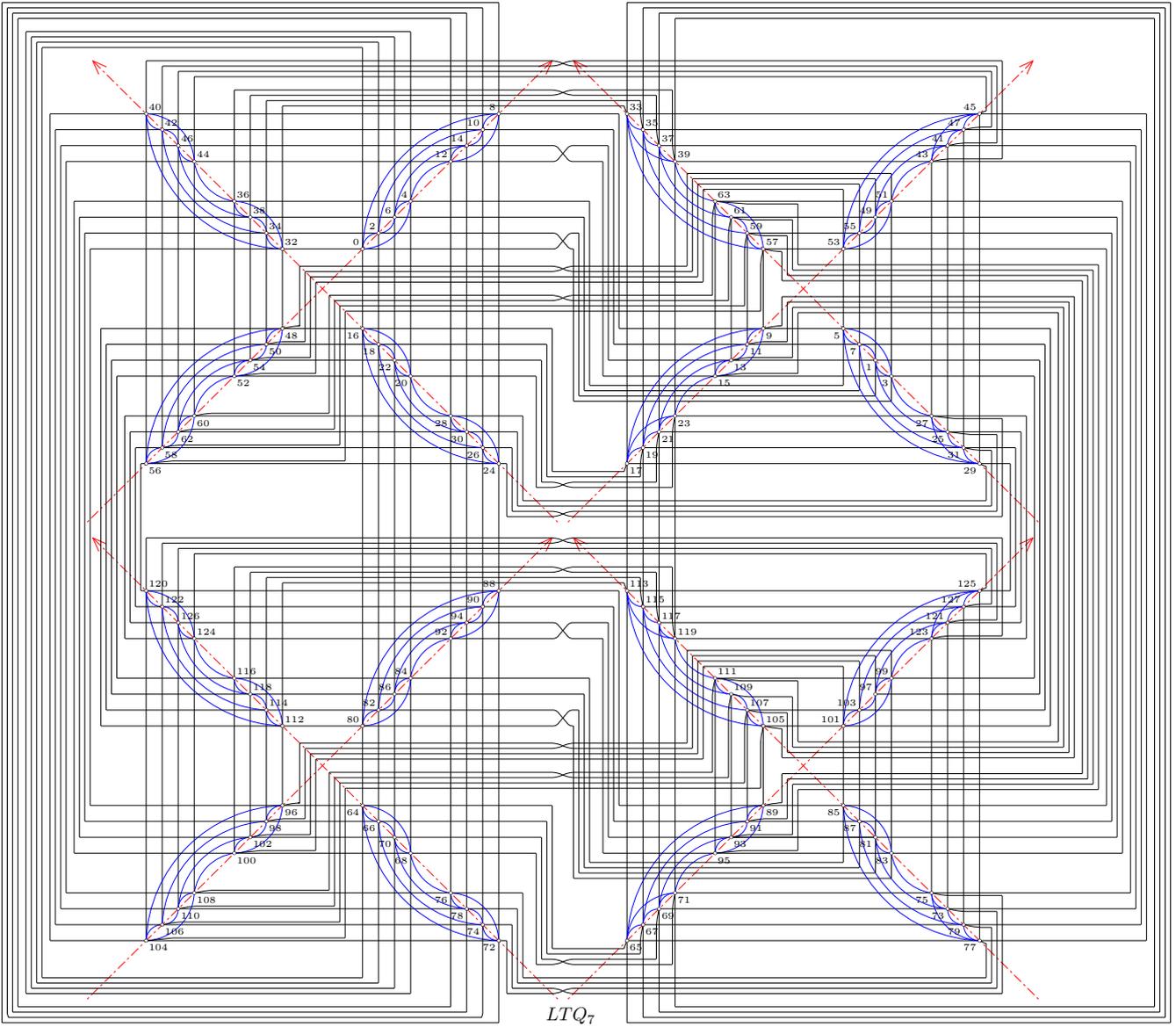}
\caption{\small{The drawing $D_7$}}
\end{figure}

\newpage

\begin{figure}[ht]
\hspace{-35pt}
\includegraphics[scale=0.625]{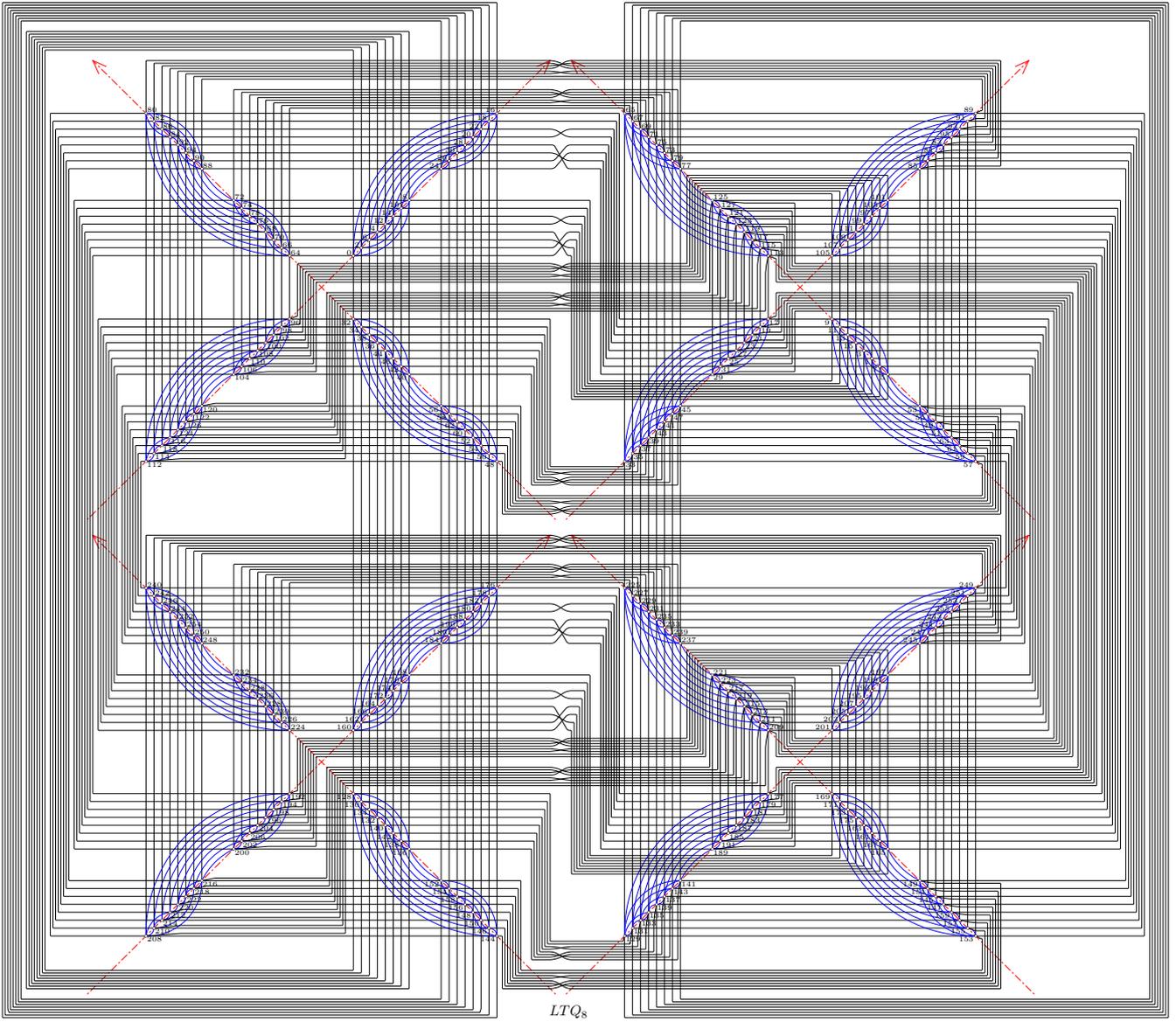}
\caption{\small{The drawing $D_8$}}
\end{figure}

\newpage

\section{Lower bound for $cr(LTQ_n)$}

\indent \indent We begin this section with the following
observation.

\begin{observation}\label{Observation the unique one in some dimension} Let $u$ be a vertex of $LTQ_n$. For any $i\in
\{1,2,\ldots, n\}$, there exists exactly one vertex $u_i \in
V(LTQ_n)$ such that $u$ and $u_i$ are adjacent with
$\lambda(u,u_i)=i$.
\end{observation}

Let $v$ be a vertex of $LTQ_n$. Let
$\tau_v:V(LTQ_n)\setminus\{v\}\rightarrow V(LTQ_n)$ be a map defined
as follows: for any vertex $u\in V(LTQ_n)\setminus\{v\}$, let
$\tau_v(u)$ be the vertex of $LTQ_n$ such that $u$ and $\tau_v(u)$
are adjacent with $\lambda(u,\tau_v(u))=\lambda(u,v)$.

It is easy to see that either $\tau_v(u)=v$ or $\lambda(u,v)+1\leq
\lambda(\tau_v(u),v)\leq n$. Hence, we can define the following.

\begin{definition} For any two vertices $u,v\in V(LTQ_n)$,
let $\mathscr{P}_{u,v}=(u_0,u_1,\ldots,u_{\ell})$ be the unique path
of $LTQ_n$ such that $u_0=u$, $u_{\ell}=v$ and $\tau_v(u_i)=u_{i+1}$
for any $i\in \{0,1,\ldots,\ell-1\}$.
\end{definition}

Note that
\begin{equation}\label{equation path dimension incresing}
\lambda(u_0,v)<\lambda(u_1,v)<\cdots <\lambda(u_{\ell-1},v).
\end{equation}

For any two vertices $v,w\in V(LTQ_n)$ and integers $1\leq t_1\leq
t_2\leq n$, let
$$D_{v}(t_1,t_2)=\{u\in V(LTQ_n)\setminus \{v\}: t_1\leq \lambda(u,v)\leq t_2\},$$
and let $$\mathcal{F}(v,w; t_1,t_2)=D_{v}(t_1,t_2)\cap \{u\in
V(LTQ_n)\setminus \{v\}:w \mbox{ is in }\mathscr{P}_{u,v} \}.$$

\begin{lemma}\label{lemma F(;)} Let $v,w$ be two vertices of $LTQ_n$, where $d=\lambda(w,v)$.
Let $k$ be an integer such that $1\leq k\leq d$. Then
$$|\mathcal{F}(v,w; k,d)|=2^{d-k}.$$
\end{lemma}

\begin{proof} By induction on $d-k$.
If $k=d$, it follows from \eqref{equation path dimension incresing}
that $\mathcal{F}(v,w; d,d)=\{w\}$, done. Hence, we assume $$k<d.$$
By \eqref{equation path dimension incresing}, we have
$\mathcal{F}(v,w;k,k)=\{u\in D_v(k,k):\tau_v(u)\in
\mathcal{F}(v,w;k+1,d)\}.$ Combining with Observation
\ref{Observation the unique one in some dimension}, we conclude that
$|\mathcal{F}(v,w;k,k)|=|\mathcal{F}(v,w;k+1,d)|$. It follows from
the induction hypothesis that
$|\mathcal{F}(v,w;k,d)|=|\mathcal{F}(v,w;k,k)|+|\mathcal{F}(v,w;k+1,d)|=2\times
2^{d-(k+1)}=2^{d-k}$. The lemma follows. \end{proof}

\indent \indent We shall introduce the lower bound method proposed
by Leighton \cite{L81}. Let $G_1=(V_1,E_1)$ and $G_2=(V_2,E_2)$ be
graphs. An embedding of $G_1$ in $G_2$ is a couple of mapping
$(\varphi,\kappa)$ satisfying
$$\varphi: V_1\rightarrow V_2$$ is an injection
$$\kappa: E_1\rightarrow \{\mbox{set of all paths in $G_2$}\},$$ such that
if $uv\in E_1$ then $\kappa(uv)$ is a path between $\varphi(u)$ and
$\varphi(v)$. For any $e\in E_2$ define
$$cg_e(\varphi,\kappa)=|\{f\in E_1:e\in \kappa(f)\}|$$
and
$$cg(\varphi,\kappa)=\max\limits_{e\in E_2}\{cg_e(\varphi,\kappa)\}.$$
The value $cg(\varphi,\kappa)$ is called congestion.

Let $2K_m$ be the complete multigraph of $m$ vertices, in which
every two vertices are joined by two parallel edges.

\begin{lemma}\label{Lemma congestion} \cite{L81} Let $(\varphi,\kappa)$ be an embedding of $G_1$ in
$G_2$ with congestion $cg(\varphi,\kappa)$. Let $\Delta(G_2)$ denote
the maximal degree of $G_2$. Then
$$cr(G_2)\geq \frac{cr(G_1)}{cg^2(\varphi,\kappa)}-\frac{|V_2|}{2}\Delta^2(G_2).$$
\end{lemma}

According to Erd\H{o}s \cite{EG73} and Kainen \cite{K72}, the
following lemmas are held.

\begin{lemma}\label{Lemma crossing of 2K} \cite{EG73} $cr(K_{2^n})\geq
\frac{2^n(2^n-1)(2^n-2)(2^n-3)}{80}$.
\end{lemma}

\begin{lemma}\label{Lemma crossing 2K}\cite{K72}
$cr(2K_{2^n})=4cr(K_{2^n})$.
\end{lemma}

Now we are in a position to show the lower bound of $cr(LTQ_n)$.

\begin{theorem}\label{Theorem Lower Bound for general n}
$cr(LTQ_n)>\frac{4^n}{20}-(n^2+1)2^{n-1}.$
\end{theorem}

\begin{proof} By Lemma \ref{Lemma congestion},
Lemma \ref{Lemma crossing of 2K} and Lemma \ref{Lemma crossing 2K},
we need only to construct an embedding $(\varphi,\kappa)$ of
$2K_{2^n}$ into $LTQ_n$ with congestion $cg(\varphi,\kappa)$ at most
$2^n$. Let $\varphi$ be an arbitrary bijection of $V(2K_{2^n})$ onto
$V(LTQ_n)$. We define the mapping $\kappa$ as follows. For any two
vertices $u$ and $v$ of $LTQ_n$, take $\mathscr{P}_{u,v}$ and
$\mathscr{P}_{v,u}$ to be the images (paths) of the two parallel
edges between $\varphi^{-1}(u)$ and $\varphi^{-1}(v)$ under
$\kappa$.

Let $e=xy$ be an arbitrary edge of $LTQ_n$, where $d=Dim(e)$. It
suffices to show
$$cg_e(\varphi,\kappa)\leq 2^n.$$

Consider first the number of paths $\mathscr{P}_{u,v}$ traversing
$x$ previous $y$, denoted $p(x,y)$. Let $V_{x,y}=\{v\in
V(LTQ_n)\setminus \{x\}:\tau_v(x)=y\}$. Note that
\begin{equation}\label{equation number of p(x,y)}
p(x,y)=\sum\limits_{v\in V_{x,y}}|\{u\in V(LTQ_n)\setminus\{v\}:x
\mbox{ is in } \mathscr{P}_{u,v} \}|.
\end{equation}
We see that $v\in V_{x,y}$ if and only if,
\begin{equation}\label{equation lambda (v,x)=d}
\lambda(v,x)=d,
\end{equation}
or equivalently, $$\theta_i(v)=\theta_i(x)\mbox{ for all }i\in
\{1,2,\ldots,d-1\}\mbox{ and }\theta_d(v)=\overline{\theta_d(x)}.$$
This implies that
\begin{equation}\label{equation |Vx,y|}
|V_{x,y}|=2^{n-d}.
\end{equation}
Combined with \eqref{equation path dimension incresing},
\eqref{equation lambda (v,x)=d} and Lemma \ref{lemma F(;)}, we have
that for any $v\in V_{x,y}$,
\begin{equation}\label{equation choices of strating point}
|\{u\in V(LTQ_n)\setminus \{v\}:x \mbox{ is in
}\mathscr{P}_{u,v}\}|= |\mathcal{F}(v,x; 1,d)|=2^{d-1}.
\end{equation}
By \eqref{equation number of p(x,y)}, \eqref{equation |Vx,y|} and
\eqref{equation choices of strating point}, we have
$$p(x,y)=2^{n-1}.$$

Similarly, the number $p(y,x)$ of paths $\mathscr{P}_{u,v}$
traversing $y$ previous $x$ is $2^{n-1}$. Therefore,
$$cg_e(\varphi,\kappa)=p(x,y)+p(y,x)=2^n.$$ This completes the proof of Theorem \ref{Theorem Lower Bound for general n}.
\end{proof}

\end{document}